\documentclass[10pt,a4paper]{article}

\usepackage[utf8]{inputenc}
\usepackage[colorlinks=True]{hyperref}
\usepackage[T1]{fontenc}
\usepackage{lmodern}    
\usepackage[a4paper]{geometry}
\usepackage{amsmath, amsfonts, amssymb}
\numberwithin{equation}{section}
\usepackage{graphicx}
\usepackage{xcolor}
\usepackage{ntheorem}
\usepackage[english]{babel}
 
\theoremstyle{break}
\newtheorem{theorem}{Theorem}

\newtheorem{lemme}[theorem]{Lemma}
\newtheorem{prop}[theorem]{Proposition}

\reversemarginpar

\renewcommand{\L}{\mathcal{L}}

\newcommand{\B}{\mathcal{B}}
\newcommand{\U}{\mathcal{U}}
\newcommand{\R}{\mathcal{R}}
\newcommand{\Ca}{C^{1+\alpha}}
\newcommand{\Cb}{C^{1+\beta}}
\title{Regularity of the spectrum for expanding maps}
\date{\today}
\author{Sedro Julien\footnote{Universit\'e Paris-Sud, B\^atiment 307, Bureau 3D3, 91400 Orsay, France. Contact: Julien.Sedro@math.u-psud.fr}}

\begin{document}
	\maketitle
	\section{Introduction and main result}
	In this short note we propose to extend one of the result in \cite{moi}
	namely differentiability of a normalized eigenfunction associated to the simple, dominating eigenvalue of the weighted transfer operator for a uniformly expanding map.
	\\More precisely, given $\B$ a (infinite dimensional) Banach space, $\U\subset\B$ an open subset, a $C^1$ family $(T_u)_{u\in\U}$ of $C^{1+\alpha}$ expanding maps on a (compact, connected, finite dimensional) Riemann manifold M, and $g:M\rightarrow \mathbb{R}$ a bounded function, the weighted transfer operator $\L_{u,g}=\L_u$, acting on $C^{1+\alpha}(M)$ by $\L_{u,g}f(y)=\sum_{x\in T_u^{-1}y}g(x)f(x)$ has the following properties :
	\begin{itemize}
		\item $\L_u$ is a bounded operator on $C^{1+\alpha}(M)$, for every $\alpha\geq-1$.
		\item If $g$ is a positive function, for every $u\in\U$, every $-1<\alpha$, it admits a spectral gap, with dominating eigenvalue $\lambda_u$, eigenfunction $\phi_u\in C^{1+\alpha}(M)$, and eigenform $\ell_u\in C^{1+\alpha}(M)'$.
	\end{itemize}
Using an implicit function approach, we established that given $0\leq\beta<\alpha$, $u\in\U\mapsto\phi_u\in C^\beta(M)$ is a differentiable map, where $\phi_u$ is chosen so that $\langle \ell_u,\phi_u\rangle=1$.
The success of this approach was largely due to the spectral gap property, which allowed us to obtain some crucial estimates. The simplicity of the dominating eigenvalue was also of major importance.
\\Here, we propose to give a generalization to the whole discrete spectrum of a quasi compact operator, by studying directly the regularity (w.r.t the parameters) of the spectral projector. This extension has a clear dynamical interest, as some dynamically meaningful quantities are connected to this spectral data (e.g, the rate of convergence to the invariant measure is given by the second largest eigenvalue of the transfer operator, see \cite[Prop 1.1]{Ba00}).
\medskip

We also want to emphasize a difference between the framework in \cite{moi} and the present work: we will restrict ourselves to transfer operator acting on \emph{little Hölder spaces} $(c^r(M))_{r\leq0}$, which are the closure of $C^\infty(M)$ with respect to the classical $C^r$ norm. This is due to our need for spectrum invariance when changing the spaces on which the transfer operator acts on (cf lemma \ref{spectralinvariance}), and the fact that smooth functions are not in general dense in the spaces of Hölder functions.
\medskip

Finally, we would like to acknowledge that the result established here is already known, and can be proved by other techniques (notably \emph{weak spectral perturbation} \`a la Keller-Gou\"ezel-Liverani, of which we give a glimpse here). The approach interest lies, in our opinion, in the method of proof: it shows that the \emph{implicit function approach} first devised in \cite{moi} can be used to study the non-dominant discrete spectrum of a quasi-compact operator, a fact which was not obvious at first sight.

We can now state the main result of this note:
\begin{theorem}\label{mainresult}
	Let M be a smooth compact, connected, boundaryless Riemann manifold of finite dimension, $\mathcal{B}$ be a Banach space, and $u_0\in\mathcal{B}$. 
	\\Let $0\leq\beta<\alpha$, and $(T_u)_{u\in\U}$ be a $C^1$ family of $C^{1+\alpha}$ expanding maps, and let $g\in C^{1+\alpha}(U\times M)$ be a (positive) function. Let $\L_{g,T_u}=\L_u$ be the weighted transfer operator, acting on the little Hölder space $c^{1+\alpha}(M)$. 
	\\Let $\lambda\not\in\sigma(\L_{u_0})$, and $\R(\lambda,u_0)=(\lambda\mathbf{1}-\L_{u_0})^{-1}$ be the resolvent operator. Then there exists a (small enough) open neighborhood $\U\subset\B$ such that
	\begin{itemize}
		\item For any $u\in\U$, the resolvent $\R(\lambda,u)\in L(c^{1+\alpha}(M))$ is well-defined. Furthermore, the map $u\in\U\mapsto\R(\lambda,u)\in L(c^{1+\alpha}(M),c^{1+\beta}(M))$ is $\gamma:=\alpha-\beta$ Hölder.
		\item  $u\in\U\mapsto\R(\lambda,u)\in L(c^{1+\beta}(M),c^\beta(M))$ is differentiable, and its differential satisfies
		\begin{equation}\label{resolventdiff}
		\partial_{u}\R(\lambda,u)=-\R(\lambda,u)\partial_u\L_u\R(\lambda,u)
		\end{equation}
	\end{itemize}
\end{theorem}

\section{Proof of the main result}
As in \cite{moi}, our strategy relies on an implicit function theorem, whose statement we recall here:

\begin{theorem}\label{myresult}
	Let $\mathcal{B}$, $X_{0}, X_{1}$ be Banach spaces such that $X_{1}\hookrightarrow X_{0}$, and the injection $j_{0}:X_{1}\rightarrow X_{0}$ is a continuous, linear map. 
	\\Let $A_{1}\subset X_{1}$ be a closed subset, and $A_{0}=j_{0}(A_{1})\subset X_{0}$. 
	\\Let $u_{0}\in\mathcal{B}$, and $U$ a neighborhood of $u_{0}$ in $\mathcal{B}$. 
	\\Consider continuous maps $F_{i}:U\times A_{i}\rightarrow A_{i}$, where $i\in\{0,1\}$, with the following property :
	\begin{equation}
	F_{0}(u,j_{0}(\phi_{1}))=j_{0}\circ F_{1}(u,\phi_{1})
	\end{equation}
	for all $u\in U$, $\phi_{1}\in A_{1}$.
	
	Moreover, we will assume that :
	
	\begin{enumerate}
		
		\item [(i)] For every $u\in U$, $F_{1}(u,.):A_{1}\rightarrow A_{1}$ admits a fixed point $\phi_{1}(u)\in A_{1}$.
		\\Furthermore, the map $u\in U\longmapsto\phi_{1}(u)\in X_{1}$ is continuous. 
		\item [(ii)]Let $\phi_0(u)=j_0(\phi_1(u))$. 
		\\For some $(u_{0},\phi_{0}(u_0))\in U\times j_{0}(A_{1})$, there exists $P_{u_{0},\phi_{0}(u_0)}^{(1)}\in \mathcal{L}(\mathcal{B},X_{0})$, $Q_{u_{0},\phi_{0}(u_0)}^{(1)}\in \mathcal{L}(j_{0}(X_{1}),X_{0})$, and a neighborhood $U'\subset U$ of $u_{0}$ such that
		\begin{equation}\label{NewTaylor}
		F_{0}(u_{0}+h,\phi_{0}+z_{0})-F_{0}(u_{0},\phi_{0})=P_{u_{0},\phi_{0}}^{(1)}.h+Q_{u_{0},\phi_{0}}^{(1)}.z_{0}+(\|h\|_{\mathcal{B}}+\|z_{0}\|_{X_{0}})\epsilon(h,z_{1})
		\end{equation}
		where $\phi_0=\phi_0(u_0)$, $h\in\mathcal{B}$ satisfies $u_{0}+h\in U'$, $z_{1}\in A_{1}$, $z_{0}=j_{0}(z_{1})\in A_{0}$, and $\epsilon(h,z_{1})\overset{X_{0}}{\underset{(h,z_{1})\rightarrow(0,0)}{\longrightarrow}} 0$.
		\item [(iii)]$Id-Q_{u_{0},\phi_{0}}\in \mathcal{L}(j_{0}(X_{1}),X_{0})$ can be extended to a bounded, invertible operator of $X_{0}$ into itself.
	\end{enumerate}
	
	Then the following holds :
	
	\begin{enumerate}
		\item [(i)']Let $\phi_{0}(u)=j_{0}(\phi_{1}(u))$. The map $u\in U\mapsto \phi_{0}(u)$ is differentiable at $u=u_{0}$ 
		\item [(ii)'] Its differential satisfies 
		\begin{equation}
		\frac{d\phi_{0}}{du}(u_{0})=(Id-Q_{u_{0},\phi_{0}(u_{0})})^{-1}P_{u_{0},\phi_{0}(u_{0})}
		\end{equation}
	\end{enumerate}
\end{theorem}
We refer to \cite{moi} for a proof of this result.
\bigskip

Let M be a Riemann manifold, $\mathcal{B}$ be a Banach space, and $\mathcal{U}\subset\mathcal{B}$ be an open subset. Let $(T_u)_{u\in\U}$ be a $C^1$ family of $C^{1+\alpha}$ expanding maps, and let $g\in C^{1+\alpha}(U\times M)$ be a (positive) function. We define the weighted transfer operator $\L_{g,T_u}=\L_u$ for every $u\in\U$, $f\in C^{1+\alpha}(M)$ by

\begin{equation}\label{transferopdef}
\L_uf(y)=\sum_{x,T_ux=y}g(u,x)f(x)
\end{equation}

It is known fact that, when acting on the space $C^r(M)$, $r>0$, each $\L_u$ has a \emph{spectral gap}, which is to say that it admits a simple, dominant eigenvalue $\lambda_u$ and that the rest of the spectrum lies on a disk of radius strictly smaller than $|\lambda_u|$ (\cite{Ba00,FanJiang,Ru89}). 
\\In \cite{moi}, we shown that for $0\leq\beta<\alpha$, the map $u\in\U\longmapsto\phi_u\in C^{\beta}(M)$ was differentiable, with $\phi_u\in C^{1+\alpha}(M)$ the eigenvector associated to $\lambda_u$, normalized by $\langle\ell_u,\phi_u\rangle=1$ (where $\ell_u$ is the eigenform associated to $\lambda_u$).
\bigskip

We want to deduce from Theorem \ref{mainresult} the regularity of the spectral projector $\Pi_u$ on the eigenspace associated with $\lambda_u$. In order to achieve that goal, our strategy consists in applying Theorem \ref{myresult}, to a map constructed from the resolvent operator $\mathcal{R}(\lambda,u)=(\lambda\mathbf{1}-\L_u)^{-1}$. 
\medskip

Let $\lambda\not\in\sigma(\L_{u_0})$ and $f\in c^{r}(M)$, $r\in\{\beta,\alpha,1+\beta,1+\alpha\}$. For every $u\in\U$, define the affine map $F:\U\times c^r(M)\mapsto c^r(M)$ by

\begin{equation}\label{mapofinterest}
F(u,\phi):=\frac{1}{\lambda}f-\frac{1}{\lambda}\L_u\phi
\end{equation}

Given our choice of $\lambda$, there is a unique $\varphi_{u_0}\in c^{r}(M)$ such that $\lambda\varphi_{u_0}-\L_{u_0}\varphi_{u_0}=f$, i.e the map $F(u_0,.)$ admits $\varphi_{u_0}=\R(\lambda,u_0)f$ as a unique fixed point.

\bigskip

We will establish Theorem \ref{mainresult} by applying Theorem \ref{myresult} to the map F.

\bigskip

The first step consists in establishing both existence of $\R(\lambda,u)f\in c^{1+\alpha}(M)$ for $u\in\U$, and (Hölder) continuity of the map $u\in\U\mapsto\R(\lambda,u)f\in L(c^{1+\alpha}(M),c^{1+\beta(M)})$:

\begin{prop}\label{fracresponseprop}
	Let M be a Riemann manifold, $\mathcal{B}$ be a Banach space, and $\mathcal{U}\subset\mathcal{B}$ be an open subset. Let $(T_u)_{u\in\U}$ be a Lipschitz family of $C^{1+\alpha}$ expanding maps, and let $g\in C^{1+\alpha}(U\times M)$ be a (positive) function. 
	\\Let $\L_{g,T_u}=\L_u$ be the weighted transfer operator defined by \eqref{transferopdef}, acting on $c^{1+\alpha}(M)$. 
	\\Let $u_0\in\U$ and $\lambda\not\in\sigma(\L_{u_0})$, and $\R(\lambda,u_0)=(\lambda\mathbf{1}-\L_{u_0})^{-1}$ be the resolvent operator. Let $0\leq\beta<\alpha$, and $\gamma:=\alpha-\beta$. 
	
	Then there is a neighborhood $u_0\in\U'\subset\U$ such that $\R(\lambda,u)\in L(c^{1+\alpha}(M))$ is well defined for every $u\in\U'$, and the map $u\in\U\longmapsto\R(\lambda,u)\in L(c^{1+\alpha}(M),c^{1+\beta}(M))$ is $\gamma$-Hölder.
\end{prop}
\paragraph{Proof:} Let $\lambda\not\in\sigma(\L_{u_0})$ and  $u,v\in\U'$. Then by Keller-Liverani Theorem, $\R(\lambda,u)$ and $\R(\lambda,v)$ are well-defined on $c^{1+\alpha}(M)$.
We have this easy variant of the famous resolvent formula :

\begin{equation}\label{resformula}
\underset{\in L(c^{1+\alpha}(M),c^{1+\beta}(M))}{\underbrace{\R(\lambda,u)-\R(\lambda,v)}}=\underset{\in L(c^{1+\beta}(M))}{\underbrace{(\lambda\mathbf{1}-\L_u)^{-1}}}\underset{\in L(c^{1+\alpha}(M),c^{1+\beta}(M))}{\underbrace{(\L_u-\L_v)}}\underset{\in L(c^{1+\alpha}(M))}{\underbrace{(\lambda\mathbf{1}-\L_v)^{-1}}}
\end{equation}

This formula allows one to reduce the problem of regularity for  $u\in\U\longmapsto\R(\lambda,u)\in L(c^{1+\alpha}(M),c^{1+\beta}(M))$ to the regularity of $u\in\U\longmapsto\L_u\in L(c^{1+\alpha}(M),c^{1+\beta}(M))$, which was established in \cite[Lemma 5, section 3.2]{moi}.
\begin{equation}\label{holdercontinuity}
\|\L_u-\L_v\|_{C^{1+\alpha},C^{1+\beta}}\leq C(u,v)\|u-v\|^{\gamma}
\end{equation}
with $C(u,v)$ a constant, depending on the Lipschitz and Hölder norms of the invert branches of the dynamic, which can be chosen uniform in $(u,v)$ if they lie in a (small enough) neighborhood of a $u_0\in\U$ (cf. the proof in\cite[Lemma 5, section 3.2]{moi}). 
\\Interpreting $\R(\lambda,v)$ as an operator on $\Ca$, $\L_u-\L_v$ as an operator from $\Ca$ to $\Cb$, and $\R(\lambda,u)$ as an operator on $\Cb$, one gets for every $f\in\Ca$, every $u,v\in\U$
\begin{align*}
\|(\R(\lambda,u)-\R(\lambda,v))f\|_{\Cb}&\leq \|\R(\lambda,u)\|_{\Cb}\|\L_u-\L_v\|_{\Ca,\Cb}\|\R(\lambda,v)f\|_{\Ca}\\
&\leq C(u,v)\|\R(\lambda,u)\|_{\Cb}\|\R(\lambda,v)\|_{\Ca}\|f\|_{\Ca}\|u-v\|^{\gamma}
\end{align*}

Under those conditions, one gets
\begin{equation}
\|(\R(\lambda,u)-\R(\lambda,v))\|_{\Ca,\Cb}\leq C(u_0)\|\R(\lambda,u)\|_{\Cb}\|\R(\lambda,v)\|_{\Ca}\|u-v\|^{\gamma}
\end{equation}
so that $u\in\U\mapsto\R(\lambda,u)\in L(c^{1+\alpha}(M),c^{1+\beta}(M))$ is (locally) Hölder, and therefore continuous at every $u_0\in\U$, if one can get \emph{uniform} (in $u$) bounds on the $\Ca,~\Cb$ norms of $\R(\lambda,u)$. The tool to establish that kind of bound is (a part) of the \emph{Keller-Liverani theorem}~\footnote{Note that we use a non-standard version of this result: here the parameter is multi (possibly infinite) dimensional, and the injection between $(X,\|.\|)$ and $(X,|.|)$ is supposed compact, which simplifies the set of assumptions}:
\begin{theorem} [\cite{KL98}]
	Let $(X,\|.\|)$ be a Banach space, endowed with $|.|$ a second norm, such that $Id:(X,\|.\|)\hookrightarrow(X,|.|)$ is compact. For every bounded (for the $\|.\|$ norm) operator $Q$, let $\|Q\|:=\sup\{|Qf|, f\in X, \|f\|\leq 1\}$. Let $(\B,\|.\|_{\B})$ be a second Banach space, $u_0\in\B$ and $\U\in\B$ be an open neighborhood of $u_0$.
\\Let $(P_u)_{u\in\U}$ be a family of bounded (for the $\|.\|$ norm) operators such that :
\begin{enumerate}
	\item There exists $C_1,~M$ positive constants (independent of $u$) such that for every $u\in\U$, $n\in\mathbb{N}$, $$|P_u^n|\leq C_1 M^n$$
	\item There exists $C_2,C_3$ positive constants, and $m\in(0,1)$, $m<M$, such that for every $u\in\U$, $n\in\mathbb{N}$, $$\|P_u^n f\|\leq C_2m^n\|f\|+C_3M^n|f|$$
	\item There exists a function $\tau:\U\times\U\rightarrow [0,\infty)$, such that $\tau(u,v)\underset{u\rightarrow v}{\rightarrow}0$ and $$\|P_u-P_{v}\|\leq\tau(u,v)$$ for $u,v\in\U$.
	\end{enumerate}

Let $\delta>0,~m<r<M$, and $V_{\delta,r}:=\{z\in\mathbb{C},|z|\leq r,~or~d(z,\sigma(P_0))\geq\delta\}$. \\Then there are positive constants $a=a(r)>0$ and $b=b(\delta,r)>0$ and a neighborhood $\U'(\delta,r)\subset\U$ such that for every $u\in\U'$, every $z\in V_{\delta,r}^c$, \begin{equation}\label{uniformbound}
\|(z-P_u)^{-1}f\|\leq a\|f\|+b|f|
\end{equation}
\end{theorem}

Using this result with $\|.\|=\|.\|_{C^{1+\theta}},~|.|=\|.\|_{C^1}$ ($\theta\in[\beta,\alpha]$) will give the necessary uniforms bounds on $\|\R(\lambda,u)\|_{\Ca},~\|\R(\lambda,v)\|_{\Cb}$.

\bigskip
 
The second step of the proof of Theorem \ref{mainresult} consists in establishing that $(u,\phi)\in\U\times c^{1+\alpha}(M) \mapsto\frac{1}{\lambda}\L_u\phi+\frac{1}{\lambda}f$ satisfies \eqref{NewTaylor}. But $\lambda,f$ are independent of $(u,\phi)$, so that it is enough to establish \eqref{NewTaylor} for $(u,\phi)\in\U\times c^{1+\alpha}(M)\mapsto\L_u\phi\in c^{\beta}(M)$. This is what was done in \cite[Lemma 5, §3.2]{moi}, in the space $C^\beta(M)$. It extends unchanged to the present setting, as the norms used are the same.
\bigskip

The last step of the proof consists in showing the invertibility (and boundedness of the inverse) for the partial differential (w.r.t $\phi$) of the map $(u,\phi)\in\U\times c^{1+\alpha}(M)\mapsto\frac{1}{\lambda}\L_u\phi+\frac{1}{\lambda}f$ in the space $c^\beta(M)$. But it is trivial to see that this partial differential is just $Id-\frac{1}{\lambda}\L_u$, so that its invertibility on $c^\beta$ boils down to the invariance of peripheral spectrum for $\L_u$ when changing spaces to $c^\beta(M)$ from $c^{1+\alpha}(M)$. This last fact is a consequence of the following result (its proof can be found in \cite[Appendix A, p.212, lemma A.3]{Ba16}).

\begin{lemme}\label{spectralinvariance}
Let $\B$ be a separated topological linear space, and let $(\B_1,\|.\|_1)$ and $(\B_2,\|.\|_2)$ be Banach spaces, continuously embedded in $\B$. Suppose there is a subspace $\B_0\subset\B_1\cap\B_2$ that is dense in both $(\B_1,\|.\|_1),~(\B_2,\|.\|_2)$. Let $\L:\B\rightarrow\B$ be a linear continuous map, preserving $\B_0,~\B_1,~\B_2$, and such that its restriction to $\B_1,~\B_2$ are bounded operators, whose essential spectral radii are both strictly smaller than some $\rho>0$. Then the eigenvalues of $\L|_{\B_1}$ and $\L|_{\B_2}$ on $\{z\in\mathbb{C},~|z|>\rho\}$ coincide. Furthermore, the corresponding generalized eigenspaces coincide and are contained in $\B_1\cap\B_2$.
\end{lemme}

This allows us to conclude that for every $f\in c^{1+\alpha}(M)$, the map $u\in\U\mapsto \R(\lambda,u)f\in c^\beta(M)$ is differentiable, i.e point-wise differentiability of the resolvent operator viewed as an element of $L(c^{1+\alpha}(M),c^\beta(M))$.
\medskip

A natural question is therefore the possibility to extend that result in the operator topology, i.e to show differentiability of $u\in\U\mapsto \R(\lambda,u)\in L(c^{1+\alpha}(M), c^\beta(M))$. Writing \eqref{NewTaylor} for the map $F$ at $(u_0,\R(\lambda,u_0)f)$, one gets
\begin{align}
&\left[\R(\lambda,u_0+h)-\R(\lambda,u_0)+\R(\lambda,u_0)P_{u_0}\R(\lambda,u_0)\right]f\\
\notag&=\epsilon(\|h\|_\B,\|\R(\lambda,u_0+h)f-\R(\lambda,u_0)f\|_{C^{1+\beta}})\left[\|h\|_\B+\|\left(\R(\lambda,u_0+h)-\R(\lambda,u_0)\right)f\|_{C^\beta} \right]
\end{align}

Looking at the error term $\epsilon(\|h\|_\B,\|\R(\lambda,u_0+h)f-\R(\lambda,u_0)f\|_{C^{1+\beta}})$, we get from \cite[(4.1.12)]{moi} that it is bounded by
\begin{equation}\label{boundonerrorterm}
M\|h\|^{1+\rho}+M'\|h\|^{1+\gamma}+C'_1\|h\|.\|z\|_{C^{1+\alpha}}+C'_2\|h\|.\|z\|_{C^\alpha}=[\|h\|+\|z\|_{C^\alpha}]\epsilon(h,z_{1+\alpha})
\end{equation}
where $z_{1+\alpha}$ is $\R(\lambda,u_0+h)f-\R(\lambda,u_0)f$ in $C^{1+\alpha}$ topology and $\epsilon(h,z_{1+\alpha})\underset{(h,z_{1+\alpha})\rightarrow 0}{\longrightarrow}0$ in $C^\beta(X)$.

In particular, it is bounded uniformly in $\|f\|_{C^{1+\alpha}}$. Therefore, we obtain
\begin{align}
\sup_{\|f\|_{C^{1+\alpha}}\leq 1}\|\left[\R(\lambda,u_0+h)-\R(\lambda,u_0)+\R(\lambda,u_0)P_{u_0}\R(\lambda,u_0)\right]f\|_{C^\beta}=\|h\|_\B\epsilon(h)
\end{align}
i.e $u\mapsto\R(\lambda,u)\in L(c^{1+\alpha}(M),c^\beta(M))$ is differentiable at $u_0\in\U$. \hfill$\square$
\bigskip

The final part of this note is to deduce regularity of the spectral projector map $u\in\U\mapsto\Pi_u\in L(c^{1+\alpha}(M),c^{\beta}(M))$.

\begin{theorem}
	Let M be a smooth compact, connected, boundary-less Riemann manifold of finite dimension, $\mathcal{B}$ be a Banach space, and $\mathcal{U}\subset\mathcal{B}$ be an open subset. Let $0\leq\beta<\alpha$, and $(T_u)_{u\in\U}$ be a $C^1$ family of $C^{2+\alpha}$ expanding maps, and let $g\in C^{1+\alpha}(U\times M)$ be a (positive) function. Let $\L_{g,T_u}=\L_u$ be the weighted transfer operator, acting on the \emph{little Hölder space} $c^{1+\alpha}(M)$. Let $\lambda_u\in\sigma(\L_u)$ be an isolated eigenvalue of finite multiplicity, and $\Pi_u$ be the associated spectral projector on the generalized eigenspace. Then we have the following :
	\begin{itemize}
		\item The map $u\in\U\mapsto\Pi_u\in L(c^{1+\alpha}(M),c^{1+\beta}(M))$ is $\gamma:=\alpha-\beta$ Hölder.
		\item $u\in\U\mapsto\Pi_u\in c^\beta(M)$ is differentiable.
	\end{itemize}
\end{theorem}
\textbf{Proof:}
Recall the formula defining a spectral projector: denoting by $\lambda$ an isolated eigenvalue of a bounded operator $\L$, and letting $\mathcal{C}$ be a circle centered at $\lambda$, small enough not to encounter the rest of the spectrum, one has  $$\Pi=\frac{1}{2i\pi}\int_{\mathcal{C}}(z-\L)^{-1}dz$$
In our setting, this reads $\Pi_u=\frac{1}{2i\pi}\int_{\mathcal{C}_u}\R(\lambda,u)d\lambda$, with $\mathcal{C}_u$ any closed curve encircling $\lambda_u$ (and no other element of $\L_u$'s spectrum). The most natural idea now is to use a "regularity under the integral" result. 
\bigskip

The first step is to pick a closed curve that is independent of $u\in\U$, at least in a neighborhood of a $u_0\in\U$. Let $u_0\in\U$, and take $u$ in a (small enough) neighborhood $U'$ of $u_0$. It follows from Keller-Liverani Theorem that whenever $u\in U'\mapsto \L_u\in L(c^{1+\alpha}(M),c^{1+\beta}(M))$ is (uniformly) continuous, then any finite set of isolated points $\lambda_1(u),..,\lambda_k(u)$ in $\L_u$'s discrete spectrum also depends (uniformly) continuously on $u\in U'$. Therefore one can pick a radius $r_{u_0}$ such that for every $u\in U'$, $\lambda_u\in D(\lambda_{u_0},\frac{r_{u_{0}}}{2})$ and there is no other element of the spectrum of $(\L_u)_{u\in U'}$ in this disk. It is then natural to pick $\mathcal{C}_{u_0}=\mathcal{C}(\lambda_{u_0},\frac{r_{u_0}}{2})$ the circle centered in $\lambda_{u_0}$ of radius $\frac{r_{u_0}}{2}$.
Hence one can write, for any $u\in U'$: 
\begin{equation}
\Pi_u=\frac{1}{2i\pi}\int_{\mathcal{C}_{u_0}}\R(\lambda,u)d\lambda
\end{equation}

The second step is to bound $\partial_u\R(\lambda,u)\in L(c^{1+\alpha}(M),c^{\beta}(M))$ independently of $u\in U'$. 
\\Up to reduce further $U'$, one has the formula $\partial_{u}\R(\lambda,u)f=-\R(\lambda,u)\partial_u\L_u\R(\lambda,u)f$ for $u\in U'$. Through a standard trick, one can show that $\partial_{u}\R(\lambda,u)f-\partial_u\R(\lambda,u_0)f$ is a sum of terms involving only products of the terms $(\R(\lambda,u)-\R(\lambda,u_0))$, $(\partial_u\L_u-\partial_u\L_{u_0})$, $\R(\lambda,u_0)$ and $\partial_u\L_{u_0}$. Thus it is enough to bound the first two above (independently of $u\in U'$).
Bounding the first expression, one gets $\|\R(\lambda,u)-\R(\lambda,u_0)\|_{\Ca,C^\beta}\leq C(u_0,\lambda)diam(U')^{1+\gamma}$, thanks to estimate (3.1.18) in \cite[p.13, lemma 4]{moi}.
\\ To bound $(\partial_u\L_u-\partial_u\L_{u_0})$ in $C^\beta$-norm, note that one has the formula:
\begin{equation}
\partial_u\L_u\phi=-\L_u[\phi\frac{dg}{g}.X_u+d\phi.X_u]
\end{equation}

with $X_u(x)=(dT_u(x))^{-1}\circ\partial_uT_u(x):\B\rightarrow L(T_{T_ux}M,T_xM)$ is linear and bounded for every $x\in M$. Therefore it is enough to establish that :
\begin{enumerate}
	\item [(1)] $\|\L_u-\L_{u_0}\|_{C^{1+\alpha},C^\beta}$ is bounded independently of $u\in U'$.
	\item [(2)] $X_u$ can be bounded independently of $u\in U'$ 
\end{enumerate}

The first item is just another avatar of estimate (3.1.18) in \cite[p.13, lemma 4]{moi}. The second naturally follows from the regularity hypothesis made on $u\in\U\mapsto T_u$, as one has immediately: $$\|X_u(x)\|\leq \sup_{u\in\U}\|DT_u^{-1}\|_{L(T_{T_ux}M,T_xM)}\|\partial_uT_u(x)\|_{L(\B,T_{T_ux}M)}$$

One then gets a bound for $\|\partial_u\R(\lambda,u)\|_{C^\beta}$, independent of $u$ in $U'$, and trivially integrable in $\lambda$ on $\mathcal{C}_{u_0}$. Therefore, one can apply the Lebesgue's theorem for differentiability under the integral, and get the announced result.
\bibliographystyle{plain}
\bibliography{/home/sedro/Dropbox/These/biblio}

\end{document}